\documentclass[12pt]{article}
\usepackage{latexsym,amssymb,amsmath}
\begin{document}
\baselineskip=20pt

\newcommand{\la}{\langle}
\newcommand{\ra}{\rangle}
\newcommand{\psp}{\vspace{0.4cm}}
\newcommand{\pse}{\vspace{0.2cm}}
\newcommand{\ptl}{\partial}
\newcommand{\dlt}{\delta}
\newcommand{\sgm}{\sigma}
\newcommand{\al}{\alpha}
\newcommand{\be}{\beta}
\newcommand{\G}{\Gamma}
\newcommand{\gm}{\gamma}
\newcommand{\vs}{\varsigma}
\newcommand{\Lmd}{\Lambda}
\newcommand{\lmd}{\lambda}
\newcommand{\td}{\tilde}
\newcommand{\vf}{\varphi}
\newcommand{\yt}{Y^{\nu}}
\newcommand{\wt}{\mbox{wt}\:}
\newcommand{\rd}{\mbox{Res}}
\newcommand{\ad}{\mbox{ad}}
\newcommand{\stl}{\stackrel}
\newcommand{\ol}{\overline}
\newcommand{\ul}{\underline}
\newcommand{\es}{\epsilon}
\newcommand{\dmd}{\diamond}
\newcommand{\clt}{\clubsuit}
\newcommand{\vt}{\vartheta}
\newcommand{\ves}{\varepsilon}
\newcommand{\dg}{\dagger}
\newcommand{\tr}{\mbox{Tr}}
\newcommand{\ga}{{\cal G}({\cal A})}
\newcommand{\hga}{\hat{\cal G}({\cal A})}
\newcommand{\Edo}{\mbox{End}\:}
\newcommand{\for}{\mbox{for}}
\newcommand{\kn}{\mbox{ker}}
\newcommand{\Dlt}{\Delta}
\newcommand{\rad}{\mbox{Rad}}
\newcommand{\rta}{\rightarrow}
\newcommand{\mbb}{\mathbb}
\newcommand{\lra}{\Longrightarrow}
\newcommand{\X}{{\cal X}}
\newcommand{\Y}{{\cal Y}}
\newcommand{\Z}{{\cal Z}}
\newcommand{\U}{{\cal U}}
\newcommand{\V}{{\cal V}}
\newcommand{\W}{{\cal W}}
\setlength{\unitlength}{3pt}

\begin{center}{\Large \bf Polynomial Representation of $F_4$ and a New }\end{center}
\begin{center}{\Large \bf  Combinatorial Identity about Twenty-Four}
\end{center}
\vspace{0.2cm}

\begin{center}{\large Xiaoping Xu}\end{center}
\begin{center}{Institute of Mathematics, Academy of Mathematics \& System Sciences}\end{center}
\begin{center}{Chinese Academy of Sciences, Beijing 100190, P.R. China}
\footnote{Research supported
 by China NSF 10871193}\end{center}

\vspace{0.6cm}

 \begin{center}{\Large\bf Abstract}\end{center}

\vspace{1cm} {\small Singular vectors of a representation of a
finite-dimensional simple Lie algebra are weight vectors in the
underlying module that are nullified by positive root vectors. In
this article, we use partial differential equations to find all the
singular vectors of the polynomial representation of the simple Lie
algebra of type $F_4$ over its basic irreducible module. As
applications, we obtain a new combinatorial identity about the
number 24 and explicit generators of invariants. Moreover,
 we show that the number of irreducible submodules contained in the
 space of homogeneous harmonic polynomials with degree $k\geq 2$ is $\geq
 [\!|k/3|\!]+[\!|(k-2)/3|\!]+2$.}

\section{Introduction}

It has been known for may years that the representation theory of
Lie algebra is closely related to combinatorial identities.
Macdonald [M] generalized the Weyl denominator identities for finite
root systems to those for infinite affine root systems, which are
now known as the Macdonald's identities. Lepowsky and Garland [LG]
gave a homologic proof of Macdonald's identities. Kac (e.g., cf
[Ka]) derived these identities from his generalization of Weyl's
character formula for the integrable representations of affine
Kac-Moody algebras, known as Weyl-Kac formula. Lepowsky and Wilson
[LW1, LW2] found a representation theoretic proof of the
Rogers-Ramanujan identities. There are a number of the other works
relating combinatorial identities to representations of Lie
algebras.

We present here a consequence  of the Macdonald's identities taken
from Kostant's work [Ko1]. Let ${\cal G}$ be a finite-dimensional
simple Lie algebra over the field $\mbb{C}$ of complex numbers.
Denote by $\Lmd^+$ the set of dominant weights of ${\cal G}$ and by
$V(\lmd)$ the finite-dimensional irreducible ${\cal G}$-module with
highest weight $\lmd$. It is known that the Casimir operator takes a
constant $c(\lmd)$ on $V(\lmd)$. Macdonald's Theorem implies that
there exists a map $\chi:\Lmd^+\rta\{-1,0,1\}$ such that the
following identity holds:
$$(\prod_{n=1}^\infty(1-q^n))^{\dim {\cal
G}}=\sum_{\lmd\in\Lmd^+}\chi(\lmd)(\dim
V(\lmd))q^{c(\lmd)}.\eqno(1.1)$$ Kostant [K2] found a connection of
the above identity with the abelian subalgebras of ${\cal G}$.

 The number 24 is
important in our life; for instance, we have 24 hours a day.
Mathematically, it is also a very special number. The minimal
length of doubly-even self-dual binary linear codes is 24. Indeed
there is a unique such code of length 24 (cf. [P]), known as the
binary Golay code (cf. [Go]). The automorphism group of this code
is a sporadic finite simple group. The minimal dimension of  even
unimodular (self-dual) integral linear lattices without elements
of square length 2 is also 24. Again there exists a unique such
lattice of dimension 24 (cf. [C1]), known as Leech's lattice (cf.
[Le]). Conway [C2] found three sporadic finite simple groups from
the automorphism group of Leech's lattice. Griess [Gr] constructed
the Monster, the largest sporadic finite simple group, as the
automorphism group of a commutative non-associative algebra
related to Leech's lattice. The Dedekind function
$$\eta(z)=q^{1/24}\prod_{n=1}^\infty(1-q^n)\;\;\mbox{with}\;\; q=e^{2\pi
zi}\eqno(1.2)$$ is a fundamental modular form of weight $1/2$ in
number theory, where $q^{1/24}$ is crucial for the modularity.
Moreover, the Ramanujan series
$$\Dlt_{24}(z)=(\eta(z))^{24}=q\prod_{n=1}^\infty(1-q^n)^{24}=\sum_{n=1}^\infty
\tau(n)q^n,\eqno(1.3)$$ where $\tau(n)$ is called the $\tau$-{\it
function} of Ramanujan. The function $\tau(n)$ is multiplicative and
has nice congruence properties such as $\tau(7m+3)\equiv
0\;(\mbox{mod}\:7)$ and $\tau(23+k)\equiv 0\;(\mbox{mod}\:23)$ when
$k$ is any quadratic non-residue of 23. The theta series of an
integral linear lattice is the generating function of counting the
numbers of lattice points on spheres. Hecke [He] proved that the
theta series of any even unimodular lattice must be a polynomial in
the Essenstein series $E_4(z)$ and the Ramanujan series
$\Dlt_{24}(z)$. The factor $(1-q^n)^{24}$ is important to
$\Dlt_{24}(z)$. Let $\lmd_r$ be the $r$th fundamental weight of
${\cal G}$. In this article, we obtain the following identity
$$(1+q)(1+q+q^2)=(1-q)^{24}\sum_{m,n,k=0}^\infty (\dim
V(m\lmd_3+(n+k)\lmd_4))q^{3m+2n+k}\eqno(1.4)$$ when ${\cal G}$ is
the simple Lie algebra of type $F_4$. In other words, the
dimensions of the modules $V(k\lmd_3+l\lmd_4)$ are linearly
correlated by the binomial coefficients of 24. Numerically,
\begin{eqnarray*}\dim V(k\lmd_3+l\lmd_4)\!\!&=&\!\!\frac{(l+1)(k+3)(k+l+4)}{ 39504568320000}(2k+l+7)(3k+l+10)
(3k+2l+11)\\ &&
\times(\prod_{r=1}^5(k+r))(\prod_{s=2}^6(k+l+s))(\prod_{q=5}^9(2k+l+q)).\hspace{2.8cm}(1.5)\end{eqnarray*}
A direct elementary proof of (1.4) seems unthinkable.

The 52-dimensional exceptional simple Lie algebra ${\cal G}^{F_4}$
of type $F_4$ can be realized as the full derivation algebra of the
unique exceptional finite-dimensional simple Jordan algebra, which
is 27-dimensional (e.g., cf. [A]). The identity element spans a
one-dimensional trivial module. The quotient space of the Jordan
algebra over the trivial module forms a 26-dimensional irreducible
${\cal G}^{F_4}$-module, which is the unique ${\cal G}^{F_4}$-module
of minimal dimension (So it is called the {\it basic module of}
${\cal G}^{F_4}$). Singular (highest-weight) vectors of a
representation of ${\cal G}^{F_4}$ are weight vectors in the
underlying module that are nullified by positive root vectors. In
this article, we use partial
 partial differential equations to
find all the singular vectors in the polynomial algebra over the
basic irreducible module of ${\cal G}^{F_4}$. Then the identity
(1.4) is a consequence of the Weyl's theorem of complete
reducibility. Another corollary of our main theorem is that the
algebra of polynomial invariants over the basic module is generated
by two explicit invariants. In addition, there is also a simple
application to harmonic analysis.

Denote by $E_{r,s}$ the square matrix with 1 as its $(r,s)$-entry
and 0 as the others. The orthogonal Lie algebra
$$o(n,\mbb{C})=\sum_{1\leq r<s\leq
n}\mbb{C}(E_{r,s}-E_{s,r}).\eqno(1.6)$$ It acts on the polynomial
algebra ${\cal A}=\mbb{C}[x_1,...,x_n]$ by
$$(E_{r,s}-E_{s,r})|_{\cal
A}=x_r\ptl_{x_s}-x_s\ptl_{x_r}.\eqno(1.7)$$ Denote by ${\cal A}_k$
the subspace of homogeneous polynomials in ${\cal A}$ with degree
$k$. When $n\geq 3$, it is well known that  the subspace of harmonic
polynomials
$${\cal H}_k=\{f\in{\cal A}_k\mid
(\ptl_{x_1}^2+\cdots+\ptl_{x_n}^2)(f)=0\}\eqno(1.8)$$ forms an
irreducible $o(n,\mbb{C})$-module. The basic module of ${\cal
G}^{F_4}$ has an invariant bilinear form. So the subspace ${\cal
H}_k^{F_4}$ of homogeneous harmonic polynomials over the basic
module with degree $k$ (in different forms) also forms a
finite-dimensional ${\cal G}^{F_4}$-submodule. According to the
Weyl's theorem of complete reducibility, it is a direct sum of
irreducible submodules. The subspace ${\cal H}_1^{F_4}={\cal A}_1$
is the basic module itself. We deduce from our main theorem that
the number of irreducible summands of ${\cal H}_k^{F_4}$ is $\geq
 [\!|k/3|\!]+[\!|(k-2)/3|\!]+2$ for $k\geq 2$.

Our idea of using partial differential equations to solve Lie
algebra problems started in our earlier works [X1] and [X2] when
we tried to find functional generators for the invariants over
curvature tensor fields and for the differential invariants of
classical groups. Later we used partial differential equations to
find the explicit formulas for all the singular vectors in the
Verma modules of $sl(n,\mbb{C})$ (cf. [X3]) and $sp(4,\mbb{C})$
(cf. [X4], where the singular vectors related to Jantzen's work
[J] for general $sp(2n,\mbb{C})$ were also explicitly given). A
few years ago, we realized that decomposing the polynomial algebra
over a finite-dimensional module of a simple Lie algebra into a
direct sum of irreducible submodules is equivalent to solving the
differential equations of flag type:
$$(d_1+f_1d_2+f_2d_3+\cdots+f_{n-1}d_n)(u)=0,\eqno(1.9)$$
where $d_1,d_2,...,d_n$ are certain commuting locally nilpotent
differential operators on the polynomial algebra
$\mbb{C}[x_1,...,x_n]$ and $f_1,...,f_{n-1}$ are polynomials
satisfying
$$d_i(f_j)=0\qquad\mbox{if}\;\;i>j.\eqno(1.10)$$
In [X5], the methods of solving such equations were given. In
particular, we found new special functions by which we are able to
explicitly give the solutions of initial value problems of a large
family of constant-coefficient linear partial differential equations
in terms of their coefficients. Recently, Luo [Lu] used our methods
to obtain explicit bases of certain infinite-dimensional
non-canonical irreducible polynomial representations for classical
simple Lie algebras.  For convenience, we will use the notion
$$\ol{i,i+j}=\{i,i+1,i+2,...,i+j\}\eqno(1.11)$$
for integer $i$ and positive integer $j$ throughout this article.

Since our proof of the main theorem heavily depends on precise
explicit representation formulas, we present in Section 2 a
construction of the basic representation of ${\cal G}^{F_4}$ from
the simple Lie algebra ${\cal G}^{E_6}$ of type $E_6$. In this way,
the reader has the whole picture of our story, and it is easier for
us to track errors. The proofs of our main theorem and its
corollaries are  given in Section 3

 \section{Basic Representation of $F_4$}

We start with the root lattice construction of the simple Lie
algebra of type $E_6$. As we all known, the Dynkin diagram of $E_6$
is as follows:

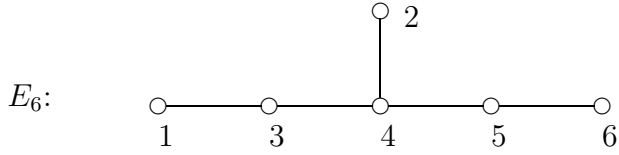
\begin{picture}(80,23)
\put(2,0){$E_6$:}\put(21,0){\circle{2}}\put(21,
-5){1}\put(22,0){\line(1,0){12}}\put(35,0){\circle{2}}\put(35,
-5){3}\put(36,0){\line(1,0){12}}\put(49,0){\circle{2}}\put(49,
-5){4}\put(49,1){\line(0,1){10}}\put(49,12){\circle{2}}\put(52,10){2}\put(50,0){\line(1,0){12}}
\put(63,0){\circle{2}}\put(63,-5){5}\put(64,0){\line(1,0){12}}\put(77,0){\circle{2}}\put(77,
-5){6}
\end{picture}
\vspace{0.7cm}

\noindent Let $\{\al_i\mid i\in\ol{1,6}\}$ be the simple positive
roots corresponding to the vertices in the diagram, and let
$\Phi_{E_6}$ be the root system of $E_6$. Set
$$Q_{E_6}=\sum_{i=1}^6\mbb{Z}\al_i,\eqno(2.1)$$ the root lattice of type
$E_6$. Denote by $(\cdot,\cdot)$ the symmetric $\mbb{Z}$-bilinear
form on $Q_{E_6}$ such that
$$\Phi_{E_6}=\{\al\in Q_{E_6}\mid (\al,\al)=2\}.\eqno(2.2)$$
 From the above Dynkin diagram of $E_6$, we have the following
automorphism of $Q_{E_6}$:
$$\sgm(\sum_{i=1}^6k_i\al_i)=k_6\al_1+k_2\al_2+k_5\al_3+k_4\al_4+k_3\al_5+k_1\al_6\eqno(2.3)$$
for $\sum_{i=1}^6k_i\al_i\in Q_{E_6}$. Define a map $F:
Q_{E_6}\times Q_{E_6}\rta \{1,-1\}$ by
$$F(\sum_{i=1}^6k_i\al_i,\sum_{j=1}^6l_j\al_j)=(-1)^{\sum_{i=1}^6k_il_i+
k_1l_3+k_4l_2+k_3l_4+k_5l_4+k_6l_5},\qquad
k_i,l_j\in\mbb{Z}.\eqno(2.4)$$ Then for $\al,\be,\gm\in Q_{E_6}$,
$$F(\al+\be,\gm)=F(\al,\gm)F(\be,\gm),\;\;F(\al,\be+\gm)=F(\al,\be)F(\al,\gm),
\eqno(2.5)$$
$$F(\al,\be)F(\be,\al)^{-1}=(-1)^{(\al,\be)},\;\;F(\al,\al)=(-1)^{(\al,\al)/2}.
\eqno(2.6)$$ In particular,
$$F(\al,\be)=-F(\be,\al)\qquad
\mbox{if}\;\;\al,\be,\al+\be\in\Phi_{E_6}.\eqno(2.7)$$ Furthermore,
$$F(\sgm(\al),\sgm(\be))=F(\al,\be)\qquad\for\;\;\al,\be\in
Q_{E_6}.\eqno(2.8)$$

Denote
$$H=\bigoplus_{i=1}^6\mbb{C}\al_i.\eqno(2.9)$$
 Then the simple  Lie algebra of type $E_6$ is
$${\cal
G}^{E_6}=H\oplus\bigoplus_{\al\in\Phi_{E_6}}\mbb{C}E_{\al}\eqno(2.10)$$
 with the Lie bracket $[\cdot,\cdot]$ determined by:
 $$[H,H]=0,\;\;[h,E_{\al}]=-[E_{\al},h]=(h,\al)E_{\al},\;\;[E_{\al},E_{-\al}]=-\al,
 \eqno(2.11)$$
 $$[E_{\al},E_{\be}]=\left\{\begin{array}{ll}0&\mbox{if}\;\al+\be\not\in\Phi_{E_6},\\
 F(\al,\be)E_{\al+\be}&\mbox{if}\;\al+\be\in\Phi_{E_6}.\end{array}\right.\eqno(2.12)$$
Moreover, we have the following automorphism $\hat\sgm$ of the Lie
algebra ${\cal G}^{E_6}$ with order 2:
$$\hat\sgm(\sum_{i=1}^6b_i\al_i)=\sum_{i=1}^6b_i\sgm(\al_i),\qquad
b_i\in\mbb{C},\eqno(2.13)$$
$$\hat\sgm(E_\al)=E_{\sgm(\al)}\qquad\for\;\;\al\in\Phi_{E_6}.\eqno(2.14)$$

The Dynkin diagram of $F_4$ is

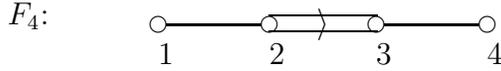
\begin{picture}(60,12)\put(2,0){$F_4$:}
\put(21,0){\circle{2}}\put(21,-5){1}\put(22,0){\line(1,0){12}}
\put(35,0){\circle{2}}\put(35,-5){2}\put(35,1.2){\line(1,0){13.6}}
\put(35,-0.8){\line(1,0){13.6}}\put(41,-1){$\ra$}\put(48.5,0){\circle{2}}\put(48.5,-5){3}\put(49.5,0)
{\line(1,0){12}}\put(62.5,0){\circle{2}}\put(62.5,-5){4}
\end{picture}
\vspace{0.6cm}

\noindent In order to make notation distinguishable, we add a bar on
the roots in the root system $\Phi_{F_4}$ of type $F_4$. In
particular, we let $\{\bar\al_1,\bar\al_2,\bar\al_3,\bar\al_4\}$ be
the simple positive roots corresponding to the above  Dynkin diagram
of $F_4$, where $\bar\al_1,\bar\al_2$ are long roots and
$\bar\al_3,\bar\al_4$ are short roots. The simple Lie algebra of
type $F_4$ is
$${\cal G}^{F_4}=\{u\in {\cal
G}^{E_6}\mid\hat\sgm(u)=u\}\eqno(2.15)$$ with Cartan subalgebra
$$H_{F_4}=\mbb{C}(\al_1+\al_6)+\mbb{C}(\al_3+\al_5)+\mbb{C}\al_4+\mbb{F}\al_2.
\eqno(2.16)$$
 Set
$$V=\{w\in {\cal
G}^{E_6}\mid\hat\sgm(w)=-w\}.\eqno(2.17)$$ Then $V$ forms the basic
26-dimensional ${\cal G}^{F_4}$-module with representation
$\ad|_{{\cal G}^{E_6}}$.

Next we want to find the explicit representation formulas of the
root vectors of ${\cal G}^{F_4}$ on $V$ in terms of differential
operators. Note the $\sgm$-invariant positive roots of $\Phi_{E_6}$
are:
$$\al_2,\;\al_4,\;\al_2+\al_4,\;\sum_{i=3}^5\al_i,\;\sum_{i=2}^5\al_i,\;
\al_1+\sum_{i=3}^6\al_i,
\;\al_4+\sum_{i=2}^5\al_i,\;\sum_{i=1}^6\al_i,\eqno(2.18)$$
$$\al_4+\sum_{i=1}^6\al_i,\;\sum_{i=1}^6\al_i+\sum_{r=3}^5\al_r,\;
\al_4+\sum_{i=1}^6\al_i+\sum_{r=3}^5\al_r,\;\al_4+\sum_{i=1}^6\al_i+\sum_{r=2}^5\al_r.
\eqno(2.19)$$ The followings are representatives of changing
positive roots modulo $\sgm$:
$$\al_1,\;\al_3,\;\al_1+\al_3,\;\al_3+\al_4,\;\al_1+\al_3+\al_4,\;
\al_2+\al_3+\al_4,\;\al_1+\al_3+\al_4+\al_5,\eqno(2.20)$$
$$\sum_{i=1}^4\al_4,\;\sum_{i=1}^5\al_i,\;\al_4+\sum_{i=1}^5\al_i,\;
\al_3+\al_4+\sum_{i=1}^5\al_i,\;\al_3+\al_4+\sum_{i=1}^6\al_i.
\eqno(2.21)$$ Set
$$x_1=E_{\al_3+\al_4+\sum_{i=1}^6\al_i}-E_{\al_5+\al_4+\sum_{i=1}^6\al_i},\qquad
x_2=E_{\al_3+\al_4+\sum_{i=1}^5\al_i}-E_{\al_4+\al_5+\sum_{i=2}^6\al_i},
\eqno(2.22)$$
$$x_3=E_{\al_4+\sum_{i=1}^5\al_i}-E_{\al_4+\sum_{i=2}^6\al_i},\qquad
x_4=E_{\sum_{i=1}^5\al_i}-E_{\sum_{i=2}^6\al_i},\eqno(2.23)$$
$$x_5=E_{\sum_{i=1}^4\al_i}-E_{\al_2+\sum_{i=4}^6\al_i},\qquad
x_6=E_{\al_1+\sum_{i=3}^5\al_i}-E_{\sum_{i=3}^6\al_i}, \eqno(2.24)$$
$$x_7=E_{\al_2+\al_3+\al_4}-E_{\al_2+\al_4+\al_5},\;\;
x_8=E_{\al_1+\al_3+\al_4}-E_{\sum_{i=4}^6\al_i},\;\;x_9=E_{\al_3+\al_4}-E_{\al_4+\al_5}
\eqno(2.25)$$
$$x_{10}=E_{\al_1+\al_3}-E_{\al_5+\al_6},\;\;x_{11}=E_{\al_3}-E_{\al_5},\;\;
x_{12}=E_{\al_1}-E_{\al_6},\eqno(2.26)$$
$$x_{13}=\al_1-\al_6,\qquad x_{14}=\al_3-\al_5,\qquad
x_{15}=E_{-\al_1}-E_{-\al_6}, \eqno(2.27)$$
$$x_{16}=E_{-\al_3}-E_{-\al_5},\;\;x_{17}=E_{-\al_1-\al_3}-E_{-\al_5-\al_6},\;\;x_{18}=E_{-\al_3-\al_4}-E_{-\al_4-\al_5},
\eqno(2.28)$$
$$x_{19}=E_{-\al_1-\al_3-\al_4}-E_{-\sum_{i=4}^6\al_i},\qquad x_{20}=
E_{-\al_2-\al_3-\al_4}-E_{-\al_2-\al_4-\al_5},\eqno(2.29)$$
$$x_{21}=E_{-\al_1-\sum_{i=3}^5\al_i}-E_{-\sum_{i=3}^6\al_i},\qquad
x_{22}=E_{-\sum_{i=1}^4\al_i}-E_{-\al_2-\sum_{i=4}^6\al_i},\eqno(2.30)$$
$$
x_{23}=E_{-\sum_{i=1}^5\al_i}-E_{-\sum_{i=2}^6\al_i},\qquad
x_{25}=E_{-\al_3-\al_4-\sum_{i=1}^5\al_i}-E_{-\al_4-\al_5-\sum_{i=2}^6\al_i},
\eqno(2.31)$$
$$x_{24}=E_{-\al_4-\sum_{i=1}^5\al_i}-E_{-\al_4-\sum_{i=2}^6\al_i},\;\;
x_{26}=E_{-\al_3-\al_4-\sum_{i=1}^6\al_i}-E_{-\al_4-\al_5-\sum_{i=1}^6\al_i}.
\eqno(2.32)$$ Then the set $\{x_i\mid i\in\ol{1,26}\}$ forms a basis
of $V$ and we treat all $x_i$ as variables for technical
convenience.

Again denote by $E_{\bar\al}$ the root vectors of ${\cal G}^{F_4}$
as follows: $\ves=\pm 1$,
$$E_{\ves\bar\al_1}=E_{\ves\al_2},\;\;E_{\ves\bar\al_2}=E_{\ves\al_4},\;\;
 E_{\ves\bar\al_3}=E_{\ves\al_3}+E_{\ves\al_5},\;\;E_{\ves\bar\al_4}=E_{\ves\al_1}+E_{\ves\al_6},
\eqno(2.33)$$
$$E_{\ves(\bar\al_1+\bar\al_2)}=E_{\ves(\al_2+\al_4)},\;\;
E_{\ves(\bar\al_2+\bar\al_3)}=E_{\ves(\al_3+\al_4)}+E_{\ves(\al_4+\al_5)}
,\eqno(2.34)$$
$$E_{\ves(\bar\al_3+\bar\al_4)}=E_{\ves(\al_1+\al_3)}+E_{\ves(\al_5+\al_6)},\;\;E_{\ves(\bar\al_1+\bar\al_2+\bar\al_3)}=
E_{\ves(\al_2+\al_3+\al_4)}+E_{\ves(\al_2+\al_4+\al_5)} ,
\eqno(2.35)$$
$$E_{\ves(\bar\al_2+\bar\al_3+\bar\al_4)}=E_{\ves(\al_1+\al_3+\al_4)}+E_{\ves(\al_4+\al_5+\al_6)},
\;\;E_{\ves(\bar\al_2+2\bar\al_3)}=E_{\ves\sum_{i=3}^5\al_i},
\eqno(2.36)$$
$$E_{\ves(\bar\al_1+\bar\al_2+2\bar\al_3)}=E_{\ves\sum_{i=2}^5\al_i},\;\;
E_{\ves(\bar\al_2+2\bar\al_3+\bar\al_4)}=E_{\ves(\al_1+\sum_{i=3}^5\al_i)}+
E_{\ves\sum_{i=3}^6\al_i},\eqno(2.37)$$
$$E_{\ves\sum_{i=1}^4\bar\al_i}=E_{\ves(\sum_{i=1}^4\al_i)}+
E_{\ves(\al_2+\sum_{i=4}^6\al_i)},\;\;E_{\ves(\bar\al_1+2\bar\al_2+2\bar\al_3)}=E_{\ves(\al_4+\sum_{i=2}^5\al_i)},
\eqno(2.38)$$
$$E_{\bar\al_3+\ves\sum_{i=1}^4\bar\al_i}=E_{\ves(\sum_{i=1}^5\al_i)}+
E_{\ves(\sum_{i=2}^6\al_i)},\;\;E_{\ves(\bar\al_2+2\bar\al_3+2\bar\al_4)}=E_{\ves(\al_1+\sum_{i=3}^6\al_i)}
\eqno(2.39)$$
$$E_{\ves(\bar\al_1+2\bar\al_2+2\bar\al_3+\bar\al_4)}=E_{\ves(\al_4+\sum_{i=1}^5\al_i)}
+E_{\ves(\al_4+\sum_{i=2}^6\al_i)},\eqno(2.40)$$
$$E_{\ves(\bar\al_1+\bar\al_2+2\bar\al_3+2\bar\al_4)}=E_{\ves\sum_{i=1}^6\al_i},\;\;
E_{\ves(\bar\al_1+2\sum_{i=2}^4\bar\al_i)}=E_{\ves(\al_4+\sum_{i=1}^6\al_i)},\eqno(2.41)$$
$$E_{\ves(\bar\al_1+2\bar\al_2+3\bar\al_3+\bar\al_4)}=E_{\ves(\al_3+\al_4+\sum_{i=1}^5\al_i)}
+E_{\ves(\al_4+\al_5+\sum_{i=2}^6\al_i)},\eqno(2.42)$$
$$E_{\ves(\bar\al_1+2\bar\al_2+3\bar\al_3+2\bar\al_4)}=E_{\ves(\al_3+\al_4+\sum_{i=1}^6\al_i)}
+E_{\ves(\al_4+\al_5+\sum_{i=1}^6\al_i)},\eqno(2.43)$$
$$E_{\ves(\bar\al_1+2\bar\al_2+4\bar\al_3+2\bar\al_4)}=E_{\ves(\sum_{i=1}^6\al_i+\sum_{r=3}^5\al_r)}
,\eqno(2.44)$$
$$E_{\ves(\bar\al_1+3\bar\al_2+4\bar\al_3+2\bar\al_4)}=E_{\ves(\al_4+\sum_{i=1}^6\al_i+\sum_{r=3}^5\al_r)}
,\eqno(2.45)$$
$$E_{\ves(2\bar\al_1+3\bar\al_2+4\bar\al_3+2\bar\al_4)}=E_{\ves(\al_4+\sum_{i=1}^6\al_i+\sum_{r=2}^5\al_r)}
.\eqno(2.46)$$ Moreover, we set
$$h_1=\al_2,\qquad h_2=\al_4,\qquad h_3=\al_3+\al_5,\qquad h_4=\al_1+\al_6.\eqno(2.47)$$

We calculate
$$E_{\bar\al_1}|_V=x_4\ptl_{x_6}+x_5\ptl_{x_8}+x_7\ptl_{x_9}-
x_{18}\ptl_{x_{20}}-x_{19}\ptl_{x_{22}}-x_{21}\ptl_{x_{23}},\eqno(2.48)$$
$$E_{\bar\al_2}|_V=x_3\ptl_{x_4}+x_8\ptl_{x_{10}}+x_9\ptl_{x_{11}}
-x_{16}\ptl_{x_{18}}-x_{17}\ptl_{x_{19}}-x_{23}\ptl_{x_{24}},\eqno(2.49)$$
\begin{eqnarray*}\hspace{1cm}E_{\bar\al_3}|_V&=&-x_2\ptl_{x_3}-x_4\ptl_{x_5}-x_6\ptl_{x_8}
+x_{10}\ptl_{x_{12}}+x_{11}(\ptl_{x_{13}}-2\ptl_{x_{14}})\\
&
&-x_{14}\ptl_{x_{16}}-x_{15}\ptl_{x_{17}}+x_{19}\ptl_{x_{21}}+x_{22}\ptl_{x_{23}}+x_{24}\ptl_{x_{25}},
\hspace{3.1cm}(2.50)\end{eqnarray*}
\begin{eqnarray*}\hspace{1cm}E_{\bar\al_4}|_V&=&-x_1\ptl_{x_2}-x_5\ptl_{x_7}-x_8\ptl_{x_9}-
x_{10}\ptl_{x_{11}}+x_{12}(\ptl_{x_{14}}-2\ptl_{x_{13}})
\\ &
&-x_{13}\ptl_{x_{15}}+x_{16}\ptl_{x_{17}}+x_{18}\ptl_{x_{19}}+x_{20}\ptl_{x_{22}}
+x_{25}\ptl_{x_{26}},\hspace{3.1cm}(2.51)\end{eqnarray*}
$$E_{\bar\al_1+\bar\al_2}|_V=-x_3\ptl_{x_6}+x_5\ptl_{x_{10}}+x_7\ptl_{x_{11}}-x_{16}\ptl_{x_{20}}
-x_{17}\ptl_{x_{22}}+x_{21}\ptl_{x_{24}},\eqno(2.52)$$
\begin{eqnarray*}\hspace{1cm}E_{\bar\al_2+\bar\al_3}|_V&=&x_2\ptl_{x_4}+x_3\ptl_{x_5}
+x_6\ptl_{10}+x_8\ptl_{x_{12}}+x_9(\ptl_{x_{13}}-2\ptl_{x_{14}})\\
&
&-x_{14}\ptl_{x_{18}}-x_{15}\ptl_{x_{19}}-x_{17}\ptl_{x_{21}}-x_{22}\ptl_{x_{24}}-x_{23}\ptl_{x_{25}},\hspace{2.5cm}
(2.53)\end{eqnarray*}
\begin{eqnarray*}\hspace{0.4cm}E_{\bar\al_3+\bar\al_4}|_V&=&-x_1
\ptl_{x_3}+x_4\ptl_{x_7}+x_6\ptl_{x_{9}}-x_{10}(\ptl_{x_{13}}+\ptl_{x_{14}})-x_{11}\ptl_{x_{15}}
\\
&&+x_{12}\ptl_{x_{16}}-(x_{13}+x_{14})\ptl_{x_{17}}-x_{18}\ptl_{x_{21}}
-x_{20}\ptl_{x_{23}}+
x_{24}\ptl_{x_{26}},\hspace{1.8cm}(2.54)\end{eqnarray*}
\begin{eqnarray*}E_{\bar\al_1+\bar\al_2+\bar\al_3}|_V&=&-x_2\ptl_{x_6}+x_3\ptl_{x_8}
+x_4\ptl_{x_{10}}+x_5\ptl_{x_{12}}+x_7(\ptl_{x_{13}}-2\ptl_{x_{14}})
-x_{14}\ptl_{x_{20}}\\ & &-x_{15}\ptl_{x_{22}}-x_{17}\ptl_{x_{23}}
-x_{19}\ptl_{x_{24}}+x_{21}\ptl_{x_{25}}, \hspace{4.7cm}(2.55)
\end{eqnarray*}
\begin{eqnarray*}E_{\bar\al_2+\bar\al_3+\bar\al_4}|_V&=&x_1\ptl_{x_4}+
x_3\ptl_{x_7}-x_6\ptl_{x_{11}}-x_8(\ptl_{x_{13}}+\ptl_{x_{14}})-x_9\ptl_{x_{15}}+x_{12}\ptl_{x_{18}}
\\ & &-(x_{13}+x_{14})\ptl_{x_{19}}+x_{16}\ptl_{x_{21}}
-x_{20}\ptl_{x_{24}}-x_{23}\ptl_{x_{26}}, \hspace{3.3cm}(2.56)
\end{eqnarray*}
$$E_{\bar\al_2+2\bar\al_3}|_V=-x_2\ptl_{x_5}+x_6\ptl_{x_{12}}+x_9\ptl_{x_{16}}-x_{11}\ptl_{x_{18}}-x_{15}\ptl_{x_{21}}
+x_{22}\ptl_{x_{25}}, \eqno(2.57)$$
$$E_{\bar\al_1+\bar\al_2+2\bar\al_3}|_V=x_2\ptl_{x_8}
+x_4\ptl_{x_{12}}+x_7\ptl_{x_{16}}-x_{11}\ptl_{x_{20}}
-x_{15}\ptl_{x_{23}}-x_{19}\ptl_{x_{25}},\eqno(2.58)$$
\begin{eqnarray*}E_{\bar\al_2+2\bar\al_3+\bar\al_4}|_V&=&-x_1\ptl_{x_5}
+x_2\ptl_{x_{7}}+x_6(\ptl_{x_{14}}-2\ptl_{x_{13}})+x_8\ptl_{x_{16}}
+x_9\ptl_{x_{17}}\\
& &-x_{10}\ptl_{x_{18}}-x_{11}\ptl_{x_{19}}-x_{13}\ptl_{x_{21}}
-x_{20}\ptl_{x_{25}}+x_{22}\ptl_{x_{26}},\hspace{2.8cm}(2.59)
\end{eqnarray*}
\begin{eqnarray*}E_{\bar\al_1+\bar\al_2+\bar\al_3+\bar\al_4}|_V&=&
-x_1\ptl_{x_6}-x_3\ptl_{x_9}-x_4\ptl_{x_{11}}-x_5(\ptl_{x_{13}}+\ptl_{x_{14}})-x_7\ptl_{x_{15}}
+x_{12}\ptl_{x_{20}}\\ &
&-(x_{13}+x_{14})\ptl_{x_{22}}+x_{16}\ptl_{x_{23}}+x_{18}\ptl_{x_{24}}
+x_{21}\ptl_{x_{26}},\hspace{2.8cm}(2.60)
\end{eqnarray*}
$$E_{\bar\al_1+2\bar\al_2+2\bar\al_3}|_V=-x_2\ptl_{x_{10}}
+x_3\ptl_{x_{12}}+x_7\ptl_{x_{18}}-x_9\ptl_{x_{20}}-x_{15}\ptl_{x_{24}}+x_{17}\ptl_{x_{25}},\eqno(2.61)$$
\begin{eqnarray*}E_{\bar\al_1+\bar\al_2+2\bar\al_3+\bar\al_4}|_V&=&x_1\ptl_{x_8}-x_2\ptl_{x_9}
+x_4(\ptl_{x_{14}}-2\ptl_{x_{13}})+x_5\ptl_{x_{16}}+x_7\ptl_{x_{17}}-x_{10}\ptl_{x_{20}}\\
&&-x_{11}\ptl_{x_{22}}-x_{13}\ptl_{x_{23}}+x_{18}\ptl_{x_{25}}-x_{19}\ptl_{x_{26}},\hspace{4cm}(2.62)
\end{eqnarray*}
$$E_{\bar\al_2+2\bar\al_3+2\bar\al_4}|_V=x_1\ptl_{x_7}
+x_6\ptl_{x_{15}}+x_8\ptl_{x_{17}}-x_{10}\ptl_{x_{19}}-x_{12}\ptl_{x_{21}}
-x_{20}\ptl_{x_{26}},\eqno(2.63)$$
\begin{eqnarray*}E_{\bar\al_1+2\bar\al_2+2\bar\al_3+\bar\al_4}|_V&=&-x_1\ptl_{x_{10}}
+x_2\ptl_{x_{11}}+x_3(\ptl_{x_{14}}-2\ptl_{x_{13}})
+x_5\ptl_{x_{18}}+x_7\ptl_{x_{19}}-x_8\ptl_{x_{20}}\\
& &-x_9\ptl_{x_{22}}-x_{13}\ptl_{x_{24}} -x_{16}\ptl_{x_{25}}
+x_{17}\ptl_{x_{26}},\hspace{4cm}(2.64)
\end{eqnarray*}
$$E_{\bar\al_1+\bar\al_2+2\bar\al_3+2\bar\al_4}|_V=
-x_1\ptl_{x_9}+x_4\ptl_{x_{15}}+x_5\ptl_{x_{17}}
-x_{10}\ptl_{x_{22}} -x_{12}\ptl_{x_{23}}
+x_{18}\ptl_{x_{26}},\eqno(2.65)$$
$$E_{\bar\al_1+2\bar\al_2+2\bar\al_3+2\bar\al_4}|_V=
x_1\ptl_{x_{11}}+x_3\ptl_{x_{15}}+x_5\ptl_{x_{19}}-x_8\ptl_{x_{22}}
-x_{12}\ptl_{x_{24}}-x_{16}\ptl_{x_{26}},\eqno(2.66)$$
\begin{eqnarray*}E_{\bar\al_1+2\bar\al_2+3\bar\al_3+\bar\al_4}|_V&=&-x_1\ptl_{x_{12}}
-x_2(\ptl_{x_{14}}+\ptl_{x_{13}})-x_3\ptl_{x_{16}}
+x_4\ptl_{x_{18}}-x_6\ptl_{x_{20}}+x_7\ptl_{x_{21}}\\
& &-x_9\ptl_{x_{23}} +x_{11}\ptl_{x_{24}}
-(x_{13}+x_{14})\ptl_{x_{25}}
+x_{15}\ptl_{x_{26}},\hspace{2.7cm}(2.67)
\end{eqnarray*}
\begin{eqnarray*}E_{\bar\al_1+2\bar\al_2+3\bar\al_3+2\bar\al_4}|_V&=&x_1(\ptl_{x_{13}}-2\ptl_{x_{14}})
+x_2\ptl_{x_{15}} -x_3\ptl_{x_{17}}+x_4\ptl_{x_{19}}
+x_5\ptl_{x_{21}}-x_6\ptl_{x_{22}}\\
& &-x_8\ptl_{x_{23}}+x_{10}\ptl_{x_{24}}-x_{12}\ptl_{x_{25}}
-x_{14}\ptl_{x_{26}},\hspace{3.9cm}(2.68)
\end{eqnarray*}
$$E_{\bar\al_1+2\bar\al_2+4\bar\al_3+2\bar\al_4}|_V=
x_1\ptl_{x_{16}}-x_2\ptl_{x_{17}} +x_4\ptl_{x_{21}}-x_6\ptl_{x_{23}}
+x_{10}\ptl_{x_{25}}-x_{11}\ptl_{x_{26}},\eqno(2.69)$$
$$E_{\bar\al_1+3\bar\al_2+4\bar\al_3+2\bar\al_4}|_V=
-x_1\ptl_{x_{18}}+x_2\ptl_{x_{19}}-x_3\ptl_{x_{21}}
+x_6\ptl_{x_{24}}-x_8\ptl_{x_{25}}+x_9\ptl_{x_{26}},\eqno(2.70)$$
$$E_{2\bar\al_1+3\bar\al_2+4\bar\al_3+2\bar\al_4}|_V=x_1\ptl_{x_{20}}
-x_2\ptl_{x_{22}} +x_3\ptl_{x_{23}} -x_4\ptl_{x_{24}}
+x_5\ptl_{x_{25}}-x_7\ptl_{x_{26}},\eqno(2.71)$$
\begin{eqnarray*}\hspace{1cm}h_1|_V&=&x_4\ptl_{x_4}+x_5\ptl_{x_5}-x_6\ptl_{x_6}+x_7\ptl_{x_7}-x_8\ptl_{x_8}-x_9\ptl_{x_9}
+x_{18}\ptl_{x_{18}}\\
&&+x_{19}\ptl_{x_{19}}-x_{20}\ptl_{x_{20}}+x_{21}\ptl_{x_{21}}-x_{22}\ptl_{x_{22}}
-x_{23}\ptl_{x_{23}},\hspace{3.4cm}(2.72)
\end{eqnarray*}
\begin{eqnarray*}\hspace{1cm}h_2|_V&=&x_3\ptl_{x_3}-x_4\ptl_{x_4}+x_8\ptl_{x_8}
+x_9\ptl_{x_9}-x_{10}\ptl_{x_{10}}-x_{11}\ptl_{x_{11}}
+x_{16}\ptl_{x_{16}}\\&
&+x_{17}\ptl_{x_{17}}-x_{18}\ptl_{x_{18}}-x_{19}\ptl_{x_{19}}+x_{23}\ptl_{x_{23}}
-x_{24}\ptl_{x_{24}},\hspace{3.4cm}(2.73)
\end{eqnarray*}
\begin{eqnarray*}\hspace{1cm}h_3|_V&=&x_2\ptl_{x_2}-x_3\ptl_{x_3}+x_4\ptl_{x_4}-x_5\ptl_{x_5}+x_6\ptl_{x_6}-x_8\ptl_{x_8}
+x_{10}\ptl_{x_{10}}\\&
&+2x_{11}\ptl_{x_{11}}-x_{12}\ptl_{x_{12}}+x_{15}\ptl_{x_{15}}-2x_{16}\ptl_{x_{16}}-x_{17}\ptl_{x_{17}}
+x_{19}\ptl_{x_{19}}
\\
&&-x_{21}\ptl_{x_{21}}+x_{22}\ptl_{x_{22}}-x_{23}\ptl_{x_{23}}+x_{24}\ptl_{x_{24}}
-x_{25}\ptl_{x_{25}},\hspace{3.4cm}(2.74)
\end{eqnarray*}
\begin{eqnarray*}\hspace{1cm}h_4|_V&=&x_1\ptl_{x_1}-x_2\ptl_{x_2}+x_5\ptl_{x_5}-x_7\ptl_{x_7}+x_8\ptl_{x_8}-x_9\ptl_{x_9}
+x_{10}\ptl_{x_{10}}\\&
&-x_{11}\ptl_{x_{11}}+2x_{12}\ptl_{x_{12}}-2x_{15}\ptl_{x_{15}}+x_{16}\ptl_{x_{16}}-x_{17}\ptl_{x_{17}}
+x_{18}\ptl_{x_{18}}
\\
&&-x_{19}\ptl_{x_{19}}+x_{20}\ptl_{x_{20}}-x_{22}\ptl_{x_{22}}+x_{25}\ptl_{x_{25}}
-x_{26}\ptl_{x_{26}}.\hspace{3.4cm}(2.75)
\end{eqnarray*}

Denote
$$\bar r=27-r\qquad\for\;\;r\in\ol{1,26}.\eqno(2.76)$$
Define an algebraic isomorphism $\tau$ on
$\mbb{C}[x_1,...,x_{26}][\ptl_{x_1},...,\ptl_{x_{26}}]$ by
$$\tau(x_r)=x_{\bar r},\;\tau(\ptl_{x_r})=\ptl_{x_{\bar
r}}\qquad\for\;\; r\in \ol{1,26}\setminus\{13,14\}\eqno(2.77)$$
and
$$\tau(x_{13})=-x_{13},\;\tau(x_{14})=-x_{14},\;\tau(\ptl_{x_{13}})=-\ptl_{x_{13}},\;\tau(\ptl_{x_{14}})=
-\ptl_{x_{14}}.\eqno(2.78)$$ Then
$$E_{-\bar\al}|_V=\tau(E_{\bar\al}|_V)\eqno(2.79)$$
for $\bar\al \in \Phi_{F_4}^+$, the set of positive roots of ${\cal
G}^{F_4}$. Thus we have given the explicit formulas for the basic
irreducible representation of $F_4$.

\section{Polynomial Representation of $F_4$}

According to (2.48)-(2.75), $V$ is the irreducible module of the
highest weight $\lmd_4$ and $x_1$ is a highest weight vector. In
this section, we want to study the ${\cal G}^{F_4}$-module ${\cal
A}=\mbb{C}[x_i\mid i\in\ol{1,26}]$ via the representation formulas
(2.48)-(2.79).

Suppose that
$$\eta_1=3\sum_{r=1}^{12}x_rx_{\bar
r}+ax_{13}^2+bx_{13}x_{14}+cx_{14}^2\eqno(3.1)$$ is an
$F_4$-invariant (cf . (2.76)), where $a,b,c$ are constants to be
determined. Then $E_{\bar\al_1}(\eta_1)=E_{\bar\al_2}(\eta_1)=0$
naturally hold. Moreover,
$$0=E_{\bar\al_3}(\eta_1)=2(a-b)x_{11}x_{13}+[b-4c-3)]x_{11}x_{14},\eqno(3.2)$$
which gives
$$a=b=4c+3, \qquad.\eqno(3.3)$$
Similarly, the constraint $E_{\bar\al_3}(\eta_1)=0$ yield
$b=c=4a+3$. So we have the quadratic invariant
$$\eta_1=3\sum_{r=1}^{12}x_rx_{\bar
r}-x_{13}^2-x_{13}x_{14}-x_{14}^2.\eqno(3.4)$$ This invariant also
gives a symmetric ${\cal G}^{F_4}$-invariant bilinear form on $V$.

By (2.72)-(2.75), we try to find a quadratic singular vector of the
form
$$\zeta_1=x_1(a_1x_{13}+a_2x_{14})+a_3x_2x_{12}+a_4x_3x_{10}+a_5x_4x_8
+a_6x_5x_6,\eqno(3.5)$$ where $a_r$ are constants. Observe
$$0=E_{\bar\al_1}(\zeta_1)=(a_5+a_6)x_4x_5\lra
a_6=-a_5.\eqno(3.5)$$ Moreover,
$$0=E_{\bar\al_2}(\zeta_1)=(a_4+a_5)x_3x_8\lra
a_5=-a_4.\eqno(3.6)$$ Note
$$0=E_{\bar\al_3}(\zeta_1)=(a_1-2a_2)x_1x_{11}+(a_3-a_4)x_2x_{10}\lra
a_1=2a_2,\;\;a_3=a_4.\eqno(3.7)$$ Furthermore,
$$0=E_{\bar\al_4}(\zeta_1)=(a_2-2a_1-a_3)x_1x_{12}\lra
2a_1=a_2-a_3.\eqno(3.8)$$ Hence we have the singular vector
$$\zeta_1=x_1(2x_{13}+x_{14})-3x_2x_{12}-3x_3x_{10}+3x_4x_8
-3x_5x_6\eqno(3.9)$$ of weight $\lmd_4$. So it generates an
irreducible module that is isomorphic to the basic module $V$.
Note
$$E_{-\bar\al_1}|_V=-x_6\ptl_{x_4}-x_8\ptl_{x_5}-x_9\ptl_{x_7}+
x_{20}\ptl_{x_{18}}+x_{22}\ptl_{x_{19}}+x_{23}\ptl_{x_{21}},\eqno(3.10)$$
$$E_{-\bar\al_2}|_V=-x_4\ptl_{x_3}-x_{10}\ptl_{x_8}-x_{11}\ptl_{x_{9}}
+x_{18}\ptl_{x_{16}}+x_{19}\ptl_{x_{17}}+x_{24}\ptl_{x_{23}},\eqno(3.11)$$
\begin{eqnarray*}\hspace{1cm}E_{-\bar\al_3}|_V&=&x_3\ptl_{x_2}+x_5\ptl_{x_4}+x_8\ptl_{x_6}
-x_{12}\ptl_{x_{10}}+x_{16}(2\ptl_{x_{14}}-\ptl_{x_{13}})\\
&&+x_{14}\ptl_{x_{11}}+x_{17}\ptl_{x_{15}}-x_{21}\ptl_{x_{19}}-x_{23}\ptl_{x_{22}}-x_{25}\ptl_{x_{24}},
\hspace{2.9cm}(3.12)\end{eqnarray*}
\begin{eqnarray*}\hspace{1cm}E_{-\bar\al_4}|_V&=&x_2\ptl_{x_1}+x_7\ptl_{x_5}+x_9\ptl_{x_8}+
x_{11}\ptl_{x_{10}}+x_{15}(2\ptl_{x_{13}}-\ptl_{x_{14}})
\\ &
&+x_{13}\ptl_{x_{12}}-x_{17}\ptl_{x_{16}}-x_{19}\ptl_{x_{18}}-x_{22}\ptl_{x_{20}}
-x_{26}\ptl_{x_{25}}\hspace{3.1cm}(3.13)\end{eqnarray*} by
(2.77)-(2.79). To get a  basis of the module generated by $\zeta_1$
compatible to $\{x_i\mid i\in\ol{1,26}\}$, we set
$$\zeta_2=E_{-\bar\al_4}(\zeta_1)=x_2(-x_{13}+x_{14})+3x_1x_{15}-3x_3x_{11}+3x_4x_9 -3x_6x_7,\eqno(3.14)$$
$$\zeta_3=E_{-\bar\al_3}(\zeta_2)=-x_3(x_{13}+2x_{14})+3x_1x_{17}+3x_2x_{16}
+3x_5x_9 -3x_7x_8,\eqno(3.15)$$
$$\zeta_4=-E_{-\bar\al_2}(\zeta_3)=-x_4(x_{13}+2x_{14})-3x_1x_{19}-3x_2x_{18} +3x_5x_{11}
-3x_7x_{10},\eqno(3.16)$$
$$\zeta_5=E_{-\bar\al_3}(\zeta_4)=x_5(-x_{13}+x_{14})+3x_1x_{21}
-3x_3x_{18}-3x_4x_{16} +3x_7x_{12},\eqno(3.17)$$
$$\zeta_6=-E_{-\bar\al_1}(\zeta_4)=-x_6(x_{13}+2x_{14})+3x_1x_{22}+3x_2x_{20}+3x_8x_{11}
-3x_9x_{10},\eqno(3.18)$$
$$\zeta_7=E_{-\bar\al_4}(\zeta_5)=x_7(2x_{13}+x_{14})+3x_2x_{21}
+3x_3x_{19}+3x_4x_{17}-3x_5x_{15}
 ,\eqno(3.19)$$
$$\zeta_8=-E_{-\bar\al_1}(\zeta_5)=x_8(-x_{13}+x_{14})
-3x_1x_{23}+3x_3x_{20}-3x_6x_{16} +3x_9x_{12},\eqno(3.20)$$
$$\zeta_9=-E_{-\bar\al_1}(\zeta_7)=x_9(2x_{13}+x_{14})-3x_2x_{23}
-3x_3x_{22} +3x_6x_{17}-3x_8x_{15},\eqno(3.21)$$
$$\zeta_{10}=-E_{-\bar\al_2}(\zeta_8)=x_{10}(-x_{13}+x_{14})
+3x_1x_{24}+3x_4x_{20}+3x_6x_{18} +3x_{11}x_{12},\eqno(3.22)$$
$$\zeta_{11}=-E_{-\bar\al_2}(\zeta_9)=x_{11}(2x_{13}+x_{14})+3x_2x_{24}
-3x_4x_{22} -3x_6x_{19}-3x_{10}x_{15},\eqno(3.23)$$
$$\zeta_{12}=-E_{-\bar\al_3}(\zeta_{10})=-x_{12}(x_{13}+2x_{14})+3x_1x_{25}
-3x_5x_{20}-3x_8x_{18} -3x_{10}x_{16},\eqno(3.24)$$
\begin{eqnarray*}\zeta_{13}=E_{-\bar\al_4}(\zeta_{12})&=&-x_{13}(x_{13}+2x_{14})
-3x_1x_{26}+3x_2x_{25}+3x_5x_{22}-3x_7x_{20}\\ &
&+3x_8x_{19}-3x_9x_{18}
+3x_{10}x_{17}-3x_{11}x_{16},\hspace{4.1cm}(3.25)
\end{eqnarray*}
\begin{eqnarray*}\zeta_{14}=E_{-\bar\al_3}(\zeta_{11})&=&x_{14}(2x_{13}+x_{14})-3x_2x_{25}+3x_3x_{24}+3x_4x_{23}
-3x_5x_{22}\\ & & +3x_6x_{21}-3x_8x_{19}-3x_{10}x_{17}+3x_{12}x_{15}
,\hspace{4.1cm}(3.26)
\end{eqnarray*}
$$\zeta_r=\tau(\zeta_{\bar r})\qquad\for\;\;r\in\ol{15,27},\eqno(3.27)$$
where $\tau$ is an algebra automorphism determined by (2.77) and
(2.78). The above construction shows that the map $x_r\mapsto
\zeta_r$ determine a module isomorphism from $V$ to the module
generated by $\zeta_1$. In particular,
$$E_{\bar\al}(x_i)=ax_j\Leftrightarrow E_{\bar\al}(\zeta_i)=a\zeta_j,\qquad
a\in\mbb{C},\;\bar\al\in\Phi_{F_4}.\eqno(3.28)$$

First
\begin{eqnarray*}\vt=(x_1\zeta_2-x_2\zeta_1)/3&=&x_1(-x_2x_{13}+x_1x_{15}-x_3x_{11}+x_4x_9
-x_6x_7)\\ & &+x_2(x_2x_{12}+x_3x_{10}-x_4x_8 +x_5x_6
)\hspace{4.1cm}(3.29)\end{eqnarray*} is a singular vector of weight
$\lmd_3$. Recall that the invariant $\eta_1$ in (3.4) define an
invariant bilinear form on $V$. Thus we have the following cubic
invariant
$$\eta_2=3\sum_{r=1}^{12}(\zeta_rx_{\bar r}+x_r\zeta_{\bar
r})-2x_{13}\zeta_{13}-x_{13}\zeta_{14}-x_{14}\zeta_{13}-2x_{14}\zeta_{14}.\eqno(3.30)$$
According to (3.9) and (3.14)-(3.27), we find
\begin{eqnarray*}\eta_2&=&9(1+\tau)[(x_2x_{12}+x_3x_{10}-x_4x_8+x_5x_6)x_{26}
+(x_3x_{11} -x_4x_9+x_6x_7)x_{25}\\ & &
+(x_7x_8-x_5x_9)x_{24}+x_{10}(x_4x_{23}+x_9x_{21})-x_{11}(x_5x_{23}+x_8x_{21}+x_{12}x_{17})
 \\ & &-x_{12}(x_7x_{22}+x_9x_{19})]
+2x_{13}^3+3x_{13}^2x_{14}-3x_{13}x_{14}^2-2x_{14}^3
+3x_1(2x_{13}+x_{14})x_{26}\\& &+3x_2(x_{13}+2x_{14})x_{25}
-3x_{13}[x_3x_{24}+x_4x_{23}+x_5x_{22}+x_6x_{21}-2x_7x_{20}+x_8x_{19}\\
&&-2x_9x_{18}+x_{10}x_{17}-2x_{11}x_{16}+x_{12}x_{15}]-3x_{14}[2x_3x_{24}+2x_4x_{23}
-x_5x_{22}\\& &+2x_6x_{21}-x_7x_{20}
-x_8x_{19}-x_9x_{18}-x_{10}x_{17} -x_{11}x_{16}+2x_{12}x_{15}]
,\hspace{2.3cm}(3.31)\end{eqnarray*} where $\tau$ is an algebra
automorphism defined in (2.77) and (2.78). Denote by $\mbb{N}$ the
set of nonnegative integers. Now we are ready to prove our main
theorem.\psp

{\bf Theorem 3.1}. {\it Any polynomial $f$ in ${\cal A}$
satisfying the system of partial differential equations
$$E_{\bar\al}(f)=0\qquad\for\;\;\bar\al\in\Phi_{F_4}^+\eqno(3.32)$$ must
be a polynomial in $x_1,\zeta_1,\vt,\eta_1,\eta_2$. In particular,
the elements
$$\{x_1^{m_1}\zeta_1^{m_2}\vt^{m_3}\eta_1^{m_4}\eta_2^{m_5}\mid
m_1,m_2,m_3,m_4,m_5\in\mbb{N}\}\eqno(3.33)$$ are linearly
independent singular vectors and any singular vector is a linear
combination of those in (3.33) with the same weight. The weight of
$x_1^{m_1}\zeta_1^{m_2}\vt^{m_3}\eta_1^{m_4}\eta_2^{m_5}$ is
$m_3\lmd_3+(m_1+m_2)\lmd_4$.}

{\it Proof.} First we note
$$x_1x_{14}=\zeta_1-2x_1x_{13}+3x_2x_{12}+3x_3x_{10}-3x_4x_8
+3x_5x_6,\eqno(3.34)$$
$$3x_1x_{15}=\zeta_2+x_2(x_{13}-x_{14})+3x_3x_{11}-3x_4x_9 +3x_6x_7,\eqno(3.35)$$
$$3x_1x_{17}=\zeta_3+x_3(x_{13}+2x_{14})-3x_2x_{16}
-3x_5x_9 +3x_7x_8,\eqno(3.36)$$
$$3x_1x_{19}=3x_5x_{11}-\zeta_4-x_4(x_{13}+2x_{14})-3x_2x_{18}
-3x_7x_{10},\eqno(3.37)$$
$$3x_1x_{21}=\zeta_5+x_5(x_{13}-x_{14})
+3x_3x_{18}+3x_4x_{16} -3x_7x_{12},\eqno(3.38)$$
$$3x_1x_{22}=\zeta_6+x_6(x_{13}+2x_{14})-3x_2x_{20}-3x_8x_{11}
+3x_9x_{10},\eqno(3.39)$$
$$3x_1x_{23}=\zeta_8+x_8(x_{13}-x_{14})
-3x_3x_{20}+3x_6x_{16} -3x_9x_{12},\eqno(3.40)$$
$$3x_1x_{24}=\zeta_{10}+x_{10}(x_{13}-x_{14})
-3x_4x_{20}-3x_6x_{18}-3x_{11}x_{12},\eqno(3.41)$$
$$3x_2x_{25}+3x_1x_{26}=\eta_1-3\sum_{r=3}^{12}x_rx_{\bar
r}+x_{13}^2+x_{13}x_{14}+x_{14}^2,\eqno(3.42)$$
\begin{eqnarray*}&&3[3(x_1x_{15}+x_3x_{11} -x_4x_9+x_6x_7)+x_2(x_{13}+2x_{14})]x_{25}
\\ & &+3[3(x_2x_{12}+x_3x_{10}-x_4x_8+x_5x_6)+x_1(2x_{13}+x_{14})]x_{26}
\\ &=&\eta_2-9x_1(x_{17}x_{24}-x_{19}x_{23}+x_{21}x_{22})-9x_2(x_{16}x_{24} -x_{18}x_{23}+x_{20}x_{21})-3x_{13}^2x_{14}
\\ & &-9(1+\tau)[(x_7x_8-x_5x_9)x_{24}+x_{10}(x_4x_{23}+x_9x_{21})-x_{11}(x_5x_{23}+x_8x_{21}+x_{12}x_{17})
\\&& -x_{12}(x_7x_{22}+x_9x_{19})]-2x_{13}^3+3x_{13}x_{14}^2+2x_{14}^3
+3x_{13}[x_3x_{24}+x_4x_{23}+x_5x_{22}+x_6x_{21}\\&
&-2x_7x_{20}+x_8x_{19}-2x_9x_{18}+x_{10}x_{17}-2x_{11}x_{16}+x_{12}x_{15}]-3x_{14}[2x_3x_{24}+2x_4x_{23}
\\& &-x_5x_{22}+2x_6x_{21}-x_7x_{20} -x_8x_{19}-x_9x_{18}-x_{10}x_{17}
-x_{11}x_{16}+2x_{12}x_{15}] \hspace{1.3cm}(3.43)\end{eqnarray*}
by (3.4), (3.9), (3.14)-(3.18), (3.20), (3.22) and (3.31). Thus
$\{x_r\mid 16,18,20\neq r\in\ol{14,26}\}$ are rational functions
in
$$\{x_r,\zeta_s,\eta_1,\eta_2\mid r\in\{\ol{1,13},16,18,20\};7,9\neq
s\in\ol{1,10}\}.\eqno(3.44)$$

Suppose that $f\in {\cal A}$ is a solution of (3.32). Write $f$ as
a rational function $f_1$ in the variables of (3.44). In the
following calculations, we will always use (3.28). By (2.71),
$$0=E_{2\bar\al_1+3\bar\al_2+4\bar\al_3+2\bar\al_4}(f_1)=x_1\ptl_{x_{20}}(f_1).
\eqno(3.45)$$ So $f_1$ is independent of $x_{20}$. Moreover,
(2.70) gives
$$0=E_{\bar\al_1+3\bar\al_2+4\bar\al_3+2\bar\al_4}(f_1)=-x_1\ptl_{x_{18}}(f_1).
\eqno(3.46)$$ Hence $f_1$ is independent of $x_{18}$. Furthermore,
(2.69) yields
$$0=E_{\bar\al_1+2\bar\al_2+4\bar\al_3+2\bar\al_4}(f_1)=x_1\ptl_{x_{16}}(f_1).
\eqno(3.47)$$ Thus $f_1$ is independent of $x_{16}$. Successively
applying (2.68), (2.67), (2.66), (2.65) and (2.63) to $f_1$, we
obtain that $f_1$ is independent of $x_{13},x_{12}, x_{11}, x_9$
and $x_7$. Therefore, $f_1$ is a rational function in
$$\{x_r,\zeta_s,\eta_1,\eta_2\mid 7,9\neq r\in\ol{1,10};7,9\neq
s\in\ol{1,10}\}.\eqno(3.48)$$

By (2.54), (2.56), (2.59), (2.60), (2.62) and (2.64),
$$0=E_{\bar\al_3+\bar\al_4}(f_1)=-x_1\ptl_{x_3}(f_1)-\zeta_1\ptl_{\zeta_3}(f_1),
\eqno(3.49)$$
$$0=E_{\bar\al_2+\bar\al_3+\bar\al_4}(f_1)=x_1\ptl_{x_4}(f_1)+\zeta_1\ptl_{\zeta_4}(f_1),
\eqno(3.50)$$
$$0=E_{\bar\al_2+2\bar\al_3+\bar\al_4}(f_1)=-x_1\ptl_{x_5}(f_1)-\zeta_1\ptl_{\zeta_5}(f_1),
\eqno(3.51)$$
$$0=E_{\bar\al_1+\bar\al_2+\bar\al_3+\bar\al_4}(f_1)=-x_1\ptl_{x_6}(f_1)-\zeta_1\ptl_{\zeta_6}(f_1),
\eqno(3.52)$$
$$0=E_{\bar\al_1+\bar\al_2+2\bar\al_3+\bar\al_4}(f_1)=x_1\ptl_{x_8}(f_1)+\zeta_1\ptl_{\zeta_8}(f_1),
\eqno(3.53)$$
$$0=E_{\bar\al_1+2\bar\al_2+2\bar\al_3+\bar\al_4}(f_1)=-x_1\ptl_{x_{10}}(f_1)-\zeta_1\ptl_{\zeta_{10}}(f_1).
\eqno(3.54)$$ Set
$$\eta_r=x_1\zeta_r-x_r\zeta_1,\;\;r\in\ol{3,6};\;\;\eta_7=x_1\zeta_8-x_8\zeta_1,
\;\;\eta_8=x_1\zeta_{10}-x_{10}\zeta_1.\eqno(3.55)$$ By the
characteristic method of solving linear partial differential
equations, we get that $f_1$ can be written as a rational function
$f_2$ in
$$\{x_r,\zeta_s,\eta_q\mid r,s=1,2;q\in\ol{1,8}\}.\eqno(3.56)$$

Next applying (2.57), (2.58) and (2.61) to $f_2$, we get
$$0=E_{\bar\al_2+2\bar\al_3}(f_2)=-(x_1\zeta_2-\zeta_1x_2)\ptl_{\eta_5}(f_2)
=-3\vt \ptl_{\eta_5}(f_2),\eqno(3.57)$$
$$0=E_{\bar\al_1+\bar\al_2+2\bar\al_3}(f_2)=3\vt\ptl_{\eta_7}(f_2),\qquad
0=E_{\bar\al_1+2\bar\al_2+2\bar\al_3}(f_2)=-3\vt\ptl_{\eta_8}(f_2)=0\eqno(3.58)$$
(cf. (3.29)). Thus $f_2$ is independent of $\eta_5$, $\eta_7$ and
$\eta_8$. Furthermore, we apply (2.50), (2.53) and (2.55) to $f_2$
and obtain
$$0=E_{\bar\al_3}(f_2)=-3\vt \ptl_{\eta_3}(f_2),\qquad 0=E_{\bar\al_2+\bar\al_3}(f_2)=3\vt
\ptl_{\eta_4}(f_2),\eqno(3.59)$$
$$0=E_{\bar\al_1+\bar\al_2+\bar\al_3}(f_2)=-3\vt
\ptl_{\eta_6}(f_2).\eqno(3.60)$$ Therefore, $f_2$ is a rational
function in $x_1,x_2,\zeta_1,\zeta_2,\eta_1,\eta_2$.  By (2.51),
$$0=E_{\bar\al_4}(f_2)=-x_1\ptl_{x_2}(f_2)-\zeta_1\ptl_{\zeta_2}(f_2).\eqno(3.61)$$
Again the characteristic method tell us that $f_2$ can be written as
a rational function $f_3$ in $x_1,\zeta_1,\vt,\eta_1,\eta_2$. Since
$f_2=f$ is a polynomial in $\{x_r\mid r\in\ol{1,26}\}$, Expressions
(3.29), (3.34), (3.42) and (3.43) imply that $f_2$ must be a
polynomial in $x_1,\zeta_1,\vt,\eta_1,\eta_2$. The other statements
follow directly.$\qquad\Box$\psp

Calculating the weights of the singular vectors in the above
theorem, we have:\psp

 {\bf Corollary 3.2}. {\it The space of polynomial ${\cal
G}^{F_4}$-invariants over its basic module is an subalgebra of
${\cal A}$ generated by $\eta_1$ and $\eta_2$}.\psp

Let $L(m_1,m_2,m_3,m_4,m_5)$ be the ${\cal G}^{F_4}$-submodule
generated by
$x_1^{m_1}\zeta_1^{m_2}\vt^{m_3}\eta_1^{m_4}\eta_2^{m_5}$. Then
$L(m_1,m_2,m_3,m_4,m_5)$ is a finite-dimensional irreducible
${\cal G}^{F_4}$-submodule with the highest weight
$m_3\lmd_3+(m_1+m_2)\lmd_4$. Let ${\cal A}_k$ be the subspace of
polynomials in ${\cal A}$ with degree $k$. Then ${\cal A}_k$ is a
finite-dimensional ${\cal G}^{F_4}$-module. By the Weyl's theorem
of complete reducibility,
$${\cal A}=\bigoplus_{k=0}^\infty {\cal A}_k=\bigoplus_{m_1,m_2,m_3,m_4,m_5=0}^\infty
L(m_1,m_2,m_3,m_4,m_5).\eqno(3.62)$$ Denote by $d(k,l)$ the
dimension of the highest weight irreducible module with the weight
$k\lmd_3+l\lmd_4$. The above equation imply the following
combinatorial identity:
$$\frac{1}{(1-t)^{26}}=\frac{1}{(1-t^2)(1-t^3)}\sum_{k_1,k_2,k_3=0}^\infty
d(k_1,k_2+k_3)t^{3k_1+2k_2+k_3}.\eqno(3.63)$$ Multiplying
$(1-t)^2$ to the above equation, we obtain a new combinatorial
identity about twenty-four:
$$\frac{1}{(1-t)^{24}}=\frac{1}{(1+t)(1+t+t^2)}\sum_{k_1,k_2,k_3=0}^\infty
d(k_1,k_2+k_3)t^{3k_1+2k_2+k_3}.\eqno(3.64)$$ Equivalently, we
have:\psp

{\bf Corollary 3.3}. {\it The dimensions $d(p,l)$ of the
irreducible module with the weights $k\lmd_3+l\lmd_4$ are linearly
correlated by the following identity:}
$$(1+t)(1+t+t^2)=(1-t)^{24}\sum_{k_1,k_2,k_3=0}^\infty
d(k_1,k_2+k_3)t^{3k_1+2k_2+k_3}.\eqno(3.65)$$ \pse

In the construction of the root system $\Phi_{F_4}$ from Euclidean
space (e.g.,cf. [Hu]),
$$(\lmd_1,\bar\al_1)=(\lmd_2,\bar\al_2)=1,\qquad
(\lmd_1,\bar\al_3)=(\lmd_2,\bar\al_4)=1/2.\eqno(3.66)$$ Recall
$\dlt=\lmd_1+\lmd_2+\lmd_3+\lmd_4$. By the dimension formula of
finite-dimensional irreducible modules of simple Lie algebras
(e.g.,cf. [Hu]),
\begin{eqnarray*}d(k,l)&=&\frac{\prod_{\bar\al\in\Phi^+_{F_4}}(\lmd+\dlt,\bar\al)}{\prod_{\bar\al\in\Phi^+_{F_4}}(\dlt,\bar\al)}
=\frac{(l+1)(k+3)(k+l+4)}{ 39504568320000}(2k+l+7)(3k+l+10)\\ &
&\times
(3k+2l+11)(\prod_{r=1}^5(k+r))(\prod_{s=2}^6(k+l+s))(\prod_{q=5}^9(2k+l+q)).
\hspace{2cm}(3.67)\end{eqnarray*}

Recall the quadratic invariants $\eta_1$ in (3.4). Dually we have
${\cal G}^{F_4}$ invariant complex Laplace operator
$$\Dlt_{F_4}=3\sum_{r=1}^{12}\ptl_{x_r}\ptl{x_{\bar
r}}-\ptl_{x_{13}}^2-\ptl_{x_{13}}\ptl_{x_{14}}-\ptl_{x_{14}}^2.\eqno(3.68)$$
Now the subspace of complex homogeneous harmonic polynomials with
degree $k$ is
$${\cal H}^{F_4}_k=\{f\in{\cal A}_k\mid
\Dlt_{F_4}(f)=0\}.\eqno(3.69)$$ Then ${\cal H}^{F_4}_1=V$. Assume
$k\geq 2$. Suppose that $k_1,k_2,m_1,m_2$ are nonnegative integers
such that
$$k_1+3k_2=m_1+3m_2+2=k.\eqno(3.70)$$
If $\Dlt_{F_4}(x_1^{k_1}\vt^{k_2})\neq 0$, then it is a singular
vector of degree $k-2$ with weight $k_2\lmd_3+k_1\lmd_4$. By
Theorem 3.1, ${\cal A}_{k-2}$ does not contain a singular vector
of such weight. A contradiction. Thus
$\Dlt_{F_4}(x_1^{k_1}\vt^{k_2})=0$. By the same reason,
$\Dlt_{F_4}(x_1^{m_1}\zeta_1\vt^{m_2})=0.$ Thus the irreducible
submodules
$$L(k_1,0,k_2,0,0),\;L(m_1,1,m_2,0,0,0)\subset {\cal
H}^{F_4}_k.\eqno(3.71)$$ This gives the following corollary:\psp

{\bf Corollary 3.4}. {\it The number of irreducible submodules
contained in the subspace ${\cal H}^{F_4}_k$ of complex homogeneous
harmonic polynomials with degree $k\geq 2$ is $\geq
 [\!|k/3|\!]+[\!|(k-2)/3|\!]+2$.}

We remark that the above conclusion implies a similar conclusion on
real harmonic polynomials for the real compact simple Lie algebras
of type $F_4$.

\vspace{0.5cm}

\noindent{\Large \bf References}

\hspace{0.5cm}

\begin{description}

\item[{[A]}] J. Adams, {\it Lectures on Exceptional Lie Groups}, The
University of Chicago Press Ltd., London, 1996.

\item[{[C1]}] J. Conway, A characterization of Leech's lattice, {\it
Invent. Math.} {\bf 69} (1969), 137-142.

\item[{[C2]}] J. Conway, Three lectures on exceptional groups, in
{\it Finite Simple Groups} (G. Higman and M. B. Powell eds.),
Chapter 7, pp. 215-247. Academy Press, London-New York, 1971.

\item[{[GL]}] H. Garland and J. Lepowsky, Lie algebra homology and
the Macdonald-Kac formulas,  {\it Invent. Math. }{\bf 34} (1976),
37-76.

\item[{[Go]}] M. Golay, Binary coding, {\it Trans. Inform. Theory}
{\bf 4} (1954), 23-28.

 \item[{[Gr]}] R. Griess, The friendly giant, {\it
Invent. Math. }{\bf 69} (1982), 1-102.

\item[{[He]}] E. Hecke, Analytische arithmetik der positiven
quadratischen formen, {\it Donske Vid Selsk.} (Mat.-Fys. Medd.) {\bf
17} (12) (1940), 1-134.

\item[{[Hu]}] J. E. Humphreys, {\it Introduction to Lie Algebras and Representation Theory},
 Springer-Verlag New York Inc., 1972.

\item[{[J]}] J. C. Jantzen, Zur charakterformel gewisser darstellungen halbeinfacher grunppen
 und Lie-algebrun, {\it Math. Z.} {\bf 140} (1974), 127-149.

\item[{[Ka]}] V. Kac, {\it Infinite-Dimensional Lie Algebras},
Birkh\"{a}ser, Boston, Inc., 1982.

\item[{[Ko1]}] B. Kostant, On Macdonald's $\eta$-function formula,
the Laplacian and generalized exponents, {\it Adv. Math.} {\bf 20}
(1976), 179-212.

\item[{[Ko2]}] B. Kostant, Powers of the Euler product and commutative
subalgebras of a complex simple Lie algebra, {\it Invent. Math.
}{\bf 158} (2004), 181-226.

\item[{[Le]}] J. Leech, Notes on sphere packings, {\it Can. J. Math.
} {\bf  19} (1967), 718-745.

\item[{[LW1]}] J. Lepowsky and R. Wilson,  A Lie theoretic
interpretation and proof of the Rogers-Ramanujan identities, {\it
Adv. Math.} {\bf 45} (1982), 21-72.

\item[{[LW2]}] J. Lepowsky and R. Wilson, The structure of standard
modules, I: universal algebras and the Rogers-Ramanujan identities,
 {\it Invent. Math. }{\bf 77} (1984), 199-290.

\item[{[Lu]}] C. Luo, Noncanonical polynomial representations of
classical Lie algebras, {\it arXiv: \\0804.0305[math.RT].}

\item[{[M]}] I. Macdonald, Affine root systems and Dedekind's
$\eta$-function, {\it Invent. Math. }{\bf 15} (1972), 91-143.

\item[{[P]}] V. Pless, On the uniqueness of the Golay codes, {\it
J. Comb. Theory} {\bf 5} (1968), 215-228.

\item[{[X1]}] X. Xu, Invariants over curvature tensor fields, {\it J. Algebra} {\bf 202}
(1998), 315-342.

\item[{[X2]}] X. Xu, Differential invariants of classical groups, {\it Duke Math. J.} {\bf 94}
(1998), 543-572.

\item[{[X3]}] X. Xu, Partial differential equations for singular
vectors of $sl(n)$, {\it arXiv: math/ \\ 0411146[math.QA].}

\item[{[X4]}] X. Xu, Differential equations for singular vectors of
$sp(2n)$, {\it Commun. Algebra} {\bf 33} (2005), 4177-4196.

\item[{[X5]}] X. Xu, Flag partial differential equations and
representations of Lie algebras, {\it Acta Appl. Math.} {\bf 102}
(2008), 149-280.

\end{description}

\end{document}